\documentclass{amsart}

\usepackage{epsfig}
\usepackage{amscd}
\usepackage{amsmath, amssymb}
\usepackage{graphics}
\usepackage{color}
\usepackage{dsfont}
\newtheorem*{main-theorem}{Main Theorem}
\newtheorem{proposition}{Proposition}[section]
\newtheorem{theorem}{Theorem}
\newtheorem*{old-theorem}{Theorem}
\newtheorem{lemma}[proposition]{Lemma}
\newtheorem{corollary}[proposition]{Corollary}

\theoremstyle{definition}

\newtheorem*{remark}{Remark}

\numberwithin{equation}{section}

\graphicspath{{/Users/hans/tex/symplectic/Dist/}{/Users/hans/tex/thesis/}}

\def\tR{\tilde{R}}
\def\tB{\tilde{B}}
\def\11{\mathbf{1}}
\def\00{\mathbf{0}}
\def\**{\mathbf{*}}
\def\II{\mathbf{I}}

\def\ZZ{{\mathbb Z}}
\def\reals{{\mathbb R}}

\def\rank{\mathrm{rank} \,}

\def\phi{\varphi}

\def\be{\begin{eqnarray*}}
\def\ee{\end{eqnarray*}}
\def\ben{\begin{eqnarray}}
\def\een{\end{eqnarray}}
\def\ker{\text{ker}}

\def\L2R{L_{\text{Rest}}^2}

\def\tC{\widetilde{C}}
\def\tR{\widetilde{R}}

\def\balpha{\bar{\alpha}}
\def\blambda{\bar{\lambda}}
\def\bi{\bar{i}}
\def\bI{\bar{I}}

\begin{document}
\title[Threshold Graphs]{The Colin de Verdi\`ere Graph Parameter for
  Threshold Graphs}
\author{Hans Christianson}
\address{Department of Mathematics, University of California, Berkeley, CA 94720 USA}
\email{hans@math.berkeley.edu}

\author{Felix Goldberg}
\address{Department of Mathematics, Technion-IIT, Technion City, Haifa 32000, ISRAEL}
\email{felixg@tx.technion.ac.il}

\begin{abstract}
We consider Schr\"odinger operators on threshold graphs and prove a
formula for the Colin de Verdi\`ere parameter in terms of the building
sequence.  We construct an optimal Colin de Verdi\`ere matrix for each connected
threshold graph $G$ of $n$ vertices.  For a large subclass of
threshold graphs we construct an alternative Colin de
Verdi\`ere matrix depending on a large parameter.  As a corollary to
this last construction, we give
estimates on the size of the non-zero eigenvalues of this matrix.
\end{abstract}
\maketitle


\section{Introduction and Statement of Results}

In this paper we consider
Schr\"odinger operators on graphs with a weight or metric on the
edges.  We present a formula for the Colin de Verdi\`ere
graph parameter for threshold graphs and show that the proof of the
formula provides an algorithm for constructing an optimal Colin de
Verdi\`ere matrix.  Finally, motivated by \cite{CR} we exploit the special structure of
the (flat) graph Laplacian for threshold graphs to provide an
alternate construction of an optimal Colin de Verdi\`ere matrix for a
large subclass of threshold graphs.

The Colin de Verdi\`{e}re number $\mu(G)$ of a graph is, roughly
speaking, a measure of the geometric complexity of a graph. It has
been introduced in \cite{CdV}. The comprehensive survey \cite{HLS} is a
reference for most of the known facts about $\mu(G)$.

Our first result is the following theorem.

\begin{theorem}
\label{main-theorem}
Let $G$ be a connected threshold graph on $n$ vertices built by
adding $i$ isolate vertices and $c$ cone vertices in some order.
\begin{enumerate}
\item\label{case:1}
If the construction sequence of $G$ is of the form ``cone, cone, $\ldots$'',
then $\mu(G)=c-1$.
\item\label{case:2}
If the construction sequence of $G$ is of the form ``cone, isolate, cone, $\ldots$'',
then $\mu(G)=c-1$.
\item\label{case:3}
If the construction sequence of $G$ is of the form ``cone, isolate, isolate, $\ldots$'',
then $\mu(G)=c$.
\end{enumerate}
Furthermore, an optimal Colin de Verdi\'{e}re matrix for $G$ can
be produced.
\end{theorem}

We say a graph $G$ is {\it threshold} if it is built inductively from
a single vertex by adding vertices one at a time
according to the following rules:

(1) Either make an edge from the new vertex to all previous vertices,
or

(2) Make no new edges.

\begin{figure}
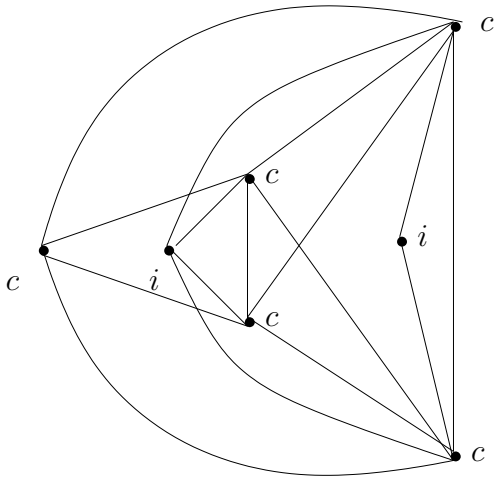

\include{fig1}
\caption{\label{fig:fig1}A threshold graph built left to right.}
\end{figure}

In case (1) we say the new vertex is a {\it cone}, and in case (2) the
new vertex is an {\it isolate}.  By convention in this work, we always
refer to the first vertex as a cone.  Hence a threshold graph might
have building sequence
\be
\text{cone, isolate, cone, cone, isolate, cone, cone}.
\ee
This graph is pictured in Figure \ref{fig:fig1}.  Clearly a threshold
graph is connected if and only if the building sequence ends with a cone.

Any two cones which appear sequentially in the building sequence are
adjacent to the same vertex set, hence are equivalent up to graph
isomorphism, and similarly for isolates which appear sequentially in
the building sequence.  Hence we may phrase the construction in terms
of {\it blocks} of cones or isolates which appear sequentially.  We
obtain a {\it block sequence}
\be
k_1, i_1, \ldots, k_m, i_m, k_{m+1}
\ee
which begins with a block of cones and ends with a block of cones under the assumption
that $G$ is connected.  The block sequence for the example in Figure
\ref{fig:fig1} is
\be
1,1,2,1,2.
\ee
Every vertex in each block has the same degree, with the last cones
added having the largest degree $d_1 = n-1$ and the last isolates having the
smallest degree $d_{2m+1} = k_{m+1}$.

There are also equivalent definitions of threshold graphs - we summarize here some of them:
\begin{theorem}[\cite{CH}]\label{thm:threshold_char}
The following statements are equivalent for a graph $G$:
\begin{enumerate}
\item
$G$ is a threshold graph.
\item
There exist weights $w_{v} \geq 0$ and a number $t$ so that for
all pairs of vertices $u \neq v$ it holds that $w_{u}+w_{v}>t$ if and only if $u$
and $v$ are adjacent.
\item
$G$ does not contain $P_{4}$, $C_{4}$ or $2K_{2}$ as an induced
subgraph.
\item\label{char:construct}
There is an assignment of  weights $w_{v} \geq 0$ to the vertices and a number $t \geq 0$
so that for a set of vertices $X$ it holds that $\sum_{v \in
X}{w_{v}} \leq t$ if and only if $X$ is independent.
\end{enumerate}
\end{theorem}

Number the vertices of a graph $G$ in order of weakly decreasing
degree, $v_1, \ldots, v_n$.  We define the graph {\it Laplacian} or
{\it incidence matrix} to be the matrix $L(G)$ given by
\be
L(G)_{ij} = \left\{ \begin{array}{l} - \# \{ \text{edges between }v_i
    \text{ and } v_j \}, \,\,\, i \neq j \\ \text{degree of } v_i, \,\,\,
    i = j. \end{array} \right.
\ee
For the example of Figure \ref{fig:fig1}, we have
\be
L(G) = \left( \begin{array}{rrrrrrr}
6 & -1 & -1 & -1 & -1 & -1 & -1 \\
-1 & 6 & -1 & -1 & -1 & -1 & -1 \\
-1 & -1 & 5 & -1 & -1 & -1 & 0 \\
-1 & -1 & -1 & 5 & -1 & -1 & 0 \\
-1 & -1 & -1 & -1 & 4 & 0 & 0 \\
-1 & -1 & -1 & -1 & 0 & 4 & 0 \\
-1 & -1 & 0 & 0 & 0 & 0 & 2
\end{array} \right).
\ee

For a less trivial example, the block sequence
\ben
\label{block-ex-2}
2,2,1,1,3,2,1
\een
results in the Laplacian
\ben
\label{lap-2}
\lefteqn{ L(G) =} \\ && \left( \begin{array}{rrrrrrrrrrrr}
11 & -1 & -1 & -1 & -1 & -1 & -1 & -1 & -1 & -1 & -1 & -1 \\

 -1 & 9 & -1 & -1 & -1 & -1 & -1 & -1 & -1 & -1 & 0 & 0
\\
 -1 & -1 & 9 & -1 & -1 & -1 & -1 & -1 & -1 & -1 & 0 & 0
\\
 -1 & -1 & -1 & 9 & -1 & -1 & -1 & -1 & -1 & -1 & 0 & 0
\\
 -1 & -1 & -1 & -1 & 8 & -1 & -1 & -1 & -1 & 0 & 0 & 0
\\
 -1 & -1 & -1 & -1 & -1 & 6 & -1 & 0 & 0 & 0 & 0 & 0 \\
 -1 & -1 & -1 & -1 & -1 & -1 & 6 & 0 & 0 & 0 & 0 & 0 \\
 -1 & -1 & -1 & -1 & -1 & 0 & 0 & 5 & 0 & 0 & 0 & 0\\
 -1 & -1 & -1 & -1 & -1 & 0 & 0 & 0 & 5 & 0 & 0 & 0 \\
 -1 & -1 & -1 & -1 & 0 & 0 & 0 & 0 & 0 & 4 & 0 & 0   \\
-1 & 0 & 0 & 0 & 0 & 0 & 0 & 0 & 0 & 0 & 1 & 0 \\
-1 & 0 & 0 & 0 & 0 & 0 & 0 & 0 & 0 & 0 & 0 &
1 \end{array} \right) . \nonumber
\een
The structure of this matrix arising from the assumption that $G$ be a
connected threshold graph is apparent in this example.  Most notably,
the matrix can be given by row blocks or column blocks corresponding
to the block sequence, and the row blocks for cones are characterized
by having no zeros before the diagonal.  We will use this in \S \ref{alt-matrix}.

Next we define {\it Colin de Verdi\`ere (CdV)} matrices, which are
edge-weighted incidence matrices plus a vertex potential satisfying
some non-degeneracy assumptions.  Specifically, we have the following
definition from \cite{HLS}:  A symmetric, real-valued $n \times n$
matrix $M$ is a CdV matrix if
\ben
 && \bullet i \neq j \implies M_{ij} <0 \text{ if } v_i \text{ and } v_j \text{
  are adjacent and } M_{ij} = 0 \text{ if not}; \label{M1} \\
 && \bullet M \text{ has exactly one negative eigenvalue of multiplicity }
1; \label{M2} \\
 && \bullet \text{there is no non-zero symmetric matrix } X \text{ satisfying }
MX=0 \label{M3} \\
&& \quad \quad \quad \text{ and } X_{ij} = 0 \text{ if } i = j \text{ or } M_{ij} \neq
0. \nonumber
\een
We think of the hypothesis \eqref{M1} as saying $M$ is a Schr\"odinger
operator on $G$, $M = L_g(G) + V(G)$,
where $L_g$ is the graph Laplacian in some Riemannian (edge-weighted)
metric $g$ and $V(G)$ is a graph potential giving weight to the
vertices.  The assumptions (\ref{M2}-\ref{M3}) correspond roughly to
saying the metric and potential are non-degenerate in some sense.

For the example of Figure \ref{fig:fig1}, the construction in the
proof of Theorem \ref{main-theorem} yields:
\ben
\label{cdv-1-felix}
\lefteqn{ M =} \nonumber \\ && \left( \begin{array}{rrrrrrr}
-0.471 &  -0.555 &  -1.02 &  -1.02 &  -0.721 &  -0.721 &  -1  \\
-0.555 &  -0.302 &  -0.474 &  -0.474 &  -0.335 &  -0.335 &  -0.129  \\
-1.02 &  -0.474 &  0 &  -0.707 &  -1 &  -1 &  0  \\
-1.02 &  -0.474 &  -0.707 &  -0.707 &  -0.5 &  -0.5 &  0  \\
-0.721 &  -0.335 &  -1 &  -0.5 &  0 &  0 &  0  \\
-0.721 &  -0.335 &  -1 &  -0.5 &  0 &  0 &  0  \\
-1 &  -0.129 &  0 &  0 &  0 &  0 &  1  \end{array} \right),
\een
while the alternative construction
in \S \ref{alt-matrix}
gives
\ben
\label{cdv-1}
M = \left( \begin{array}{rrrrrrr}
-a & -a & -b & -b & -c & -1/4 & -1/2 \\
-a & -a & -b & -b & -c & -1/4 & -1/2 \\
-b & -b & -b & -b & -c & -1/4 & 0 \\
-b & -b & -b & -b & -c & -1/4 & 0 \\
-c & -c & -c & -c & -c & 0 & 0 \\
-1/4 & -1/4 & -1/4 & -1/4 & 0 & 1/4 & 0 \\
-1/2 & -1/2 & 0 & 0 & 0 & 0 & 1/2 \end{array} \right),
\een
where $a>0$ is a sufficiently large parameter, $b = a +1/4$, and $c =
a + 3/4$.  For the block sequence \eqref{block-ex-2}, the alternative
construction from \S \ref{alt-matrix} yields
\ben
\lefteqn{M =} \nonumber \\ && -\left( \begin{array}{rrrrrrrrrrrr}
a & b & b & b & c & d & d & \frac{1}{5} & \frac{1}{5} & \frac{1}{4} & 1 & 1 \\
b & b & b & b & c & d & d & \frac{1}{5} & \frac{1}{5} & \frac{1}{4} & 0 & 0 \\
b & b & b & b & c & d & d & \frac{1}{5} & \frac{1}{5} & \frac{1}{4} & 0 & 0 \\
b & b & b & b & c & d & d & \frac{1}{5} & \frac{1}{5} & \frac{1}{4} & 0 & 0 \\
c & c & c & c & c & d & d & \frac{1}{5} & \frac{1}{5} & 0 & 0 & 0 \\
d & d & d & d & d & d & d & 0 & 0 & 0 & 0 & 0 \\
d & d & d & d & d & d & d & 0 & 0 & 0 & 0 & 0 \\
\frac{1}{5} & \frac{1}{5} & \frac{1}{5} & \frac{1}{5} & \frac{1}{5} & 0 & 0 & -\frac{1}{5} & 0 & 0 & 0 & 0 \\
\frac{1}{5} & \frac{1}{5} & \frac{1}{5} & \frac{1}{5} & \frac{1}{5} & 0 & 0 & 0 & -\frac{1}{5} & 0 & 0 & 0 \\
\frac{1}{4} & \frac{1}{4} & \frac{1}{4} & \frac{1}{4} & 0 & 0 & 0 & 0 & 0 & -\frac{1}{4} & 0 & 0 \\
1 & 0 & 0 & 0 & 0 & 0 & 0 & 0 & 0 & 0 & -1 & 0 \\
1 & 0 & 0 & 0 & 0 & 0 & 0 & 0 & 0 & 0 & 0 & -1
\end{array} \right), \label{cdv-2}
\een
where $a>0$ is a sufficiently large parameter, $b = a + 2$, $c = a +
9/4$, and $d = a + 53/20$.
That (\ref{cdv-1-felix}-\ref{cdv-2}) are CdV matrices and how they were
constructed will follow from the proof of Theorem \ref{main-theorem}
and \S \ref{alt-matrix}.

The {\it Colin de Verdi\`ere graph parameter} $\mu(G)$ is defined to be the
largest co-rank of a CdV matrix, that is the dimension of the largest
nullspace among all CdV matrices associated to $G$.  In \cite{CdV},
Colin de Verdi\`ere proved the following theorem.
\begin{theorem}\label{thm:cdv-char}
The Colin de Verdi\`ere graph parameter $\mu(G) \leq 1$ if and only if $G$ is a disjoint union of paths.

The Colin de Verdi\`ere graph parameter $\mu(G) \leq 2$ if and only if $G$ is outerplanar.

The Colin de Verdi\`ere graph parameter $\mu(G) \leq 3$ if and only if $G$ is planar.
\end{theorem}
From \cite{RST} and \cite{LS} we have the additional characterization
given by the following theorem.
\begin{theorem}
The Colin de Verdi\`ere graph parameter $\mu(G) \leq 4$ if and only if $G$ is linklessly embeddable in
$\reals^3$.
\end{theorem}
See also \cite{HLS} for a summary of these and other results.

The rank of the matrix \eqref{cdv-1} is $3$, hence $\mu(G) \geq 4$ in
this case.  As $G$ is linklessly embeddable in $\reals^3$, we have
$\mu(G) = 4$ for this graph.  The rank of the matrix \eqref{cdv-2} is
$6$, hence $\mu(G) \geq 6$ for this graph.

{\bf Acknowledgements.}  The first author would like to thank the Hill
Opportunity Fund for providing support under which he attended a
conference ({\it Festival Colin de Verdi\`ere: Semiclassical, Riemannian, and Combinatorial aspects of
Spectral Theory}) where he was introduced to the Colin de Verdi\`ere
parameter.  He would also like to thank Eran Nevo for insightful comments. The second author would like to thank Professor Abraham Berman for his kind advice and encouragement.

\section{Proof of Theorem \ref{main-theorem}}
\label{proof-section}
In this section we prove Theorem \ref{main-theorem} and construct an
optimal Colin de Verdi\`ere matrix for each threshold graph.

The following result is our main tool in the proof of
Theorem \ref{main-theorem}.
\begin{theorem}[\cite{CdV},\cite{HLS}]\label{theorem:remove}
Let $v$ be a vertex of $G$. Then $\mu(G) \leq \mu(G-v)+1$. If 
$v$ is connected to all other vertices and $G-v$ is not
$\overline{K_{2}}$ or empty, then equality holds.
\end{theorem}

We also need a construction of the optimal matrix whose existence
is asserted in the second part of Theorem \ref{theorem:remove}.
This has been given explicitly in \cite{HLS} for the case when
$G-v$ is connected. It is possible to extend this construction for
an arbitrary $G-v$ but we shall only require the case where all
connected components but one of $G-v$ are isolated vertices. Our
proof naturally follows that of \cite{HLS} with the requisite
addition.

\begin{lemma}\label{theorem:suspension_mat}
Let $G$ be a graph on $n$ vertices and let $v$ be a vertex of $G$.
Suppose that $v$ is connected to all other vertices and that $G-v$
is not $\overline{K_{2}}$ or empty. Without loss of generality
suppose that the first $k$ vertices of $G-v$ induce a connected
component $C$ of $G-v$ and that the other $n-k-1$ vertices of
$G-v$ are isolated. (Possibly, $n-k-1=0$).

Now let $M^{'}$ be a Colin de Verdi\'{e}re matrix for $C$ with
negative eigenvalue $\lambda_{1}$ . Let $z$ be a unit negative
eigenvector of $M^{'}$ corresponding to $\lambda_{1}$. Also,
without loss of generality we may assume that allf the square submatrices
associated to the $n-k-1$ isolates, 
$M_{k+1}^{'}=M_{k+2}^{'}=\ldots=M_{n-1}^{'}=1$.

Now let $\theta=\sqrt{1-(n-k-1)\lambda_{1}}$ and let
    \begin{displaymath}
    M=
        \left(
          \begin{array}{ccc}
            M^{'} & 0_{k,n-k-1} & \theta z \\
            0_{n-k-1,k} & I_{n-k-1} & -\mathbf{1}_{n-k-1} \\
            \theta z^{T} & -\mathbf{1}_{n-k-1}^{T} & \lambda_{1}^{-1} \\
          \end{array}
        \right).
    \end{displaymath}
Then $M$ is a Colin de Verdi\'{e}re matrix for $G$ and
$corank(M)=corank(M^{'})+1$.
\end{lemma}
\begin{proof}
$M$ is obviously a discrete Schr\"{o}dinger operator of $G$. Also
it is easy to verify that $M$ possesses the Strong Arnold property
since $M^{'}$ does so.

Now let $x \in \ker{M^{'}}$. The vectors of the form $(x,0,0)^{T}$
all belong to $\ker{M}$. It is not hard to verify that so does the
vector $(\theta z,\lambda \mathbf{1}^{T}_{n-k-1},\lambda)^{T}$.
Therefore $corank(M) \geq  corank(M^{'})+1$ and by the first part
of Theorem \ref{theorem:suspension_mat} the equality of coranks
follows.

The uniqueness of the negative eigenvalue of $M$ follows via
interlacing from the uniqueness of the negative eigenvalue of
$M^{'}$ and $corank(M)=corank(M^{'})+1$.
\end{proof}

We shall also need the following lemmas:
\begin{lemma}[\cite{HLS}]\label{theorem:components}
If $G$ has at least one edge then $\mu(G)=\max_{K}{\mu(K)}$, where
$K$ runs over all connected components of $G$.
\end{lemma}

\begin{lemma}\label{lemma:star}
For $q \leq 2$, $\mu(K_{1,q})=1$ and for $q \geq 3$,
$\mu(K_{1,q})=2$.

An optimal Colin de Verdi\'{e}re matrix for
$\mu(K_{1,1})$ is
 $$\left( \begin {array}{cc} -1&-1\\\noalign{\medskip}-1&-1\end {array} \right). $$

For $2 \leq q \leq 3$, an optimal Colin de Verdi\'{e}re matrix for
$\mu(K_{1,q})$ is
       $$ M=\left(
            \begin{array}{cc}
            \begin{array}{cc}0\end{array}  &
                                \begin{array}{cc}-\mathbf{1}\end{array}   \\
            \begin{array}{cc}-\mathbf{1}^{T}\end{array}  &
                                \begin{array}{cc}\mathbf{0}\end{array}
            \end{array}
        \right). $$

For $q \geq 4$, an optimal Colin de Verdi\'{e}re matrix for
$\mu(K_{1,q})$ is
      $$  M=\left(
            \begin{array}{ccc}
            \begin{array}{ccc}0\end{array}  &
                                \begin{array}{ccc}-\mathbf{1}\end{array} &
                                \begin{array}{ccc}-\mathbf{1}\end{array}  \\

            \begin{array}{ccc}-\mathbf{1}^{T}\end{array}  &
                                \begin{array}{ccc}I_{q-3}\end{array} &
                                \begin{array}{ccc}\mathbf{0}\end{array} \\

            \begin{array}{ccc}-\mathbf{1}^{T}\end{array}  &
                                \begin{array}{ccc}\mathbf{0}\end{array} &
                                \begin{array}{ccc}\mathbf{0}\end{array} \\
            \end{array}
        \right). $$

\end{lemma}
\begin{proof}
The values of $\mu$ for the stars can be read off Theorem
\ref{thm:cdv-char}. As for the matrices, they can be seen to have the
Strong Arnold Property by direct verification.
\end{proof}

\begin{proof}[Proof of Theorem \ref{main-theorem}]
Suppose that the second cone vertex in the construction sequence
has been added at stage $k$ (recall that by our convention the
sequence always starts with a cone). Therefore before stage $k$
the graph had been edgeless whereas after stage $k$ it is in fact
$K_{1,k-1}$. An optimal Colin de Verdi\'{e}re matrix for it is
given by Lemma \ref{lemma:star}.

We can trace the effect every further stage of the execution of
the construction sequence has. When we add an
isolate Lemma \ref{theorem:components} implies that we do not
alter the Colin de Verdi\'{e}re number. If $M$ is the matrix we
had constructed so far, $M \oplus (1)$ is an optimal Colin de
Verdi\'{e}re matrix after the addition of an isolate.

When we add a cone Theorem \ref{theorem:remove} implies
that we do increase the Colin de Verdi\'{e}re number by 1. An optimal Colin de
Verdi\'{e}re matrix for the graph after the addition of a cone can be now obtained
by the construction of Lemma \ref{theorem:suspension_mat}.

It remains to observe that after stage $k$ we have $c-2$ cone
additions left to do and thus the Colin de Verdi\'{e}re number of $G$
is $\mu(K_{1,k-1})+c-2$.
\end{proof}

\section{An alternative Colin de Verdi\`ere matrix}
\label{alt-matrix}
In this section we construct an alternative optimal Colin de Verdi\`ere matrix
for all but case $(3)$ in Theorem \ref{main-theorem}.
There are three
main steps, constructing a real-symmetric matrix $M$ satisfying \eqref{M1}
and proving it has
the appropriate co-rank, proving $M$ satisfies \eqref{M2}, and proving
$M$ satisfies \eqref{M3}.

\noindent {\bf Construction of $M$.}

We adopt the following labeling conventions.  By
$\11_{k \times m}$ we denote the $k \times m$ matrix of $1$s, and by
$\00_{k \times m}$ the $k \times m$ matrix of $0$s.  By $\II_{m \times
  m}$ we denote the $m \times m$ identity matrix.

Let $G$ be a graph with degree and block sequence as in the statement
of Theorem~\ref{main-theorem}.  We will construct
symmetric, real valued, $n \times n$ matrices in blocks of rows and
columns
corresponding to the block sequence of construction for a connected
threshold graph.  For
such a matrix, let $r_1, r_2, \ldots r_n$ denote the individual rows,
and $R_1, R_2, \ldots, R_{2m+1}$ denote the row blocks.  Here $R_j$
has $k_{m+2-j}$ rows for $1 \leq j \leq m+1$ corresponding to the
blocks of cones and $i_{j-m-1}$ rows for $m+2 \leq j \leq 2m+1$
corresponding to the blocks of isolates.  Let $c_1 = r_1^T,
c_2=r_2^T, \ldots c_n=r_n^T$ denote the individual columns, and
$C_1=R_1^T, C_2=R_2^T, \ldots, C_{2m+1} = R_{2m+1}^T$ denote the
column blocks.  For convenience, we will also write the expression
$r_j + R_k$ to mean ``row $j$ + $\sum$(rows in block $R_k$)'' whenever unambiguous.

We construct a family of CdV matrices for $G$.  Let $\alpha_1 ,
\ldots, \alpha_{2m+1} >0$ be a set of parameters to be fixed later in
the proof.  For our construction of $M$, first take
\be
R_1 & = & \bigg( - \alpha_1 \11_{k_{m+1} \times k_{m+1}}, - \alpha_2
  \11_{k_{m+1} \times k_m}, \ldots, -\alpha_{m+1} \11_{k_{m+1} \times
      k_1}; \\
&& \quad \quad  - \alpha_{m+2} \11_{k_{m+1} \times i_1}, \ldots, -
      \alpha_{2m+1} \11_{k_{m+1} \times i_m} \bigg),
\ee
and
\be
R_{2m+1} & = & \bigg( -\alpha_{2m+1} \11_{i_m \times k_{m+1}},
\00_{i_m \times k_m}, \ldots, \00_{i_m \times k_1} ; \\
&& \quad \quad \00_{i_m \times i_1} ,\ldots, \00_{i_m \times i_{m-1}},
\alpha_{2m+1} \II_{i_m \times i_m} \bigg).
\ee
The idea here is that $R_1$ has only one independent row, and adding
$R_{2m+1}$ to $r_1$ will kill the last block in $r_1$, which can then
be used to kill $R_2$.  More
precisely,
\be
r_1 + R_{2m+1} & = & \bigg( (- \alpha_1 - i_m \alpha_{2m+1}) \11_{1 \times
  k_{m+1}}, - \alpha_2 \11_{1 \times k_m}, \ldots - \alpha_{m+1}
\11_{1 \times k_1};  \\
&& \quad \quad - \alpha_{m+2} \11_{1 \times i_1}, \ldots, -
\alpha_{2m} \11_{1 \times i_{m-1}}, \00_{1 \times i_m} \bigg),
\ee
which is equal to $r_j$ for $k_{m+1}+1 \leq j \leq k_{m+1} + k_m$ as
long as
\ben
\label{a2-eq}
\alpha_2 = \alpha_1 + i_m \alpha_{2m+1}.
\een
Hence with this choice of $\alpha_2$, $R_2$ is dependent.  Similarly,
for $2 \leq j \leq m$, we can arrange for each row of $R_{j}$ to equal
\be
r_1 + R_{2m+1} + \ldots + R_{2m+3-j},
\ee
provided for $2 \leq j \leq m+1$,
\be
R_j & = & \bigg( - \alpha_j \11_{k_{m+2-j} \times k_{m+1}}, - \alpha_j
  \11_{k_{m+2-j} \times k_m}, \ldots, -\alpha_{j} \11_{k_{m+2-j} \times
      k_{m+2-j}}, \\
&& \quad \quad -\alpha_{j+1} \11_{k_{m+2-j}\times k_{m+3-j}}, \ldots,
-\alpha_{m} \11_{k_{m+2-j} \times k_{2}}, -\alpha_{m+1} \11_{k_{m+2-j}
\times k_1}; \\
&& \quad \quad  - \alpha_{m+2} \11_{k_{m+2-j} \times i_1}, \ldots, -
      \alpha_{2m+2-j} \11_{k_{m+2-j} \times i_{m+1-j}}, \\
&& \quad \quad  \00_{k_{m+2-j}
        \times i_{m+2-j}}, \ldots, \00_{k_{m+2-j} \times i_m} \bigg),
\ee
for $m+2 \leq j \leq 2m+1$
\be
R_j & = & \bigg( -\alpha_{j} \11_{i_{j-m-1} \times k_{m+1}},
- \alpha_j \11_{i_{j-m-1} \times k_m}, \ldots, -\alpha_j
\11_{i_{j-m-1} \times k_{j-m}} , \\
&& \quad \quad \00_{i_{j-m-1} \times k_{j-m+1}}, \ldots,
\00_{i_{j-m-1} \times k_1}; \\
&& \quad \quad \00_{i_{j-m-1} \times i_1} ,\ldots, \00_{i_{j-m-1} \times i_{j-m-2}},
\alpha_{j} \II_{i_{j-m-1} \times i_{j-m-1}},\\
&& \quad \quad \00_{i_{j-m-1} \times i_{j-m}}, \ldots, \00_{i_{j-m-1}
  \times i_m}  \bigg)
\ee
and for $2 \leq j \leq m+1$,
\ben
\label{aj-eq}
\alpha_{j} = \alpha_1 + i_m \alpha_{2m+1} + i_{m-1} \alpha_{2m} +
\ldots + i_{m+2-j} \alpha_{2m+3-j}.
\een
That is, $\alpha_2, \ldots, \alpha_{m+1}$ depend on $\alpha_1,
\alpha_{m+2}, \ldots, \alpha_{2m+1}$ and the space of available
parameters has dimension $m+1$.

We have shown
\be
\rank M = 1 + \sum_{j=1}^m i_j,
\ee
and $M$ satisfies \eqref{M1}.

{\bf Proof of property \eqref{M2}.}

For \eqref{M2} we will prove a result on the
structure of the characteristic polynomial.  We will employ
multi-index notation: Let $\epsilon \in \ZZ^l$,
\be
\epsilon = ( \epsilon_1, \ldots, \epsilon_l).
\ee
For a vector $x \in \reals^l$, by $x^\epsilon$, we mean
\be
x^\epsilon = x_1^{\epsilon_1} x_2^{\epsilon_2} \cdots
x_l^{\epsilon_l}.
\ee
By $|\epsilon|$, we mean
\be
|\epsilon| = \epsilon_1 + \ldots + \epsilon_l,
\ee
and for $\epsilon, \epsilon' \in \ZZ^{l}$, by $\epsilon \leq
\epsilon'$, we mean
\be
\epsilon_j \leq \epsilon_j'
\ee
for each $1 \leq j \leq l$, and define also $\epsilon - \epsilon'$ and
$\epsilon + \epsilon'$ componentwise.

We also use the following labeling convention when unambiguous: Let
\be
i & = & \sum_{j=1}^m i_j, \\
k & = & \sum_{j=1}^{m+1} k_j, \text{ and} \\
\alpha &=& \alpha_{m+1} = \alpha_1 + \sum_{j=1}^m i_j \alpha_{m+1+j},
\ee
and define corresponding vectors
\be
\bi & = & (i_1, \ldots, i_m), \text{ and}\\
\balpha & = & ( \alpha_{m+2}, \ldots , \alpha_{2m+1} ).
\ee
Observe $\balpha$ does not have an entry for $\alpha_1$, and has the
same length as $\bi$.  We define also the vector
\be
\blambda & = & ( \underbrace{ \lambda, \ldots, \lambda}_{m}),
\ee
so that for a multi-index $\epsilon \in \ZZ^m$, we have
\be
(\balpha - \blambda)^\epsilon = (\alpha_{m+2}- \lambda)^{\epsilon_1}
\cdots (\alpha_{2m+1} - \lambda)^{ \epsilon_m}.
\ee

To aid in computation, we introduce another set of parameters.  For $1
\leq j \leq m$, define
\be
\beta_j = k_{m+1} + k_m + \ldots + k_{j+1},
\ee
so that $\beta_j < \beta_{j-1}$ and
\be
\beta_{j-p}-\beta_{j} = k_{j} + \ldots + k_{j-p+1}.
\ee
We choose
\ben
\label{choose-alpha}
\alpha_{m+1+j} = \frac{1}{\beta_j},
\een
leaving $\alpha_1$ free.

Let $M$ be the matrix constructed above.  We will calculate $\det (M -
\lambda \II)$ by first using a similarity transformation to produce
rows of zeros in $M$, and then
using properties of the $\det$ function.

\begin{proposition}
\label{poly-prop}
Suppose $M$ is a real symmetric matrix constructed according to the
algoritheorem above with this choice of the $\alpha_j$.  Then if
$\alpha_1>0$ is chosen sufficiently large, there exist positive
constants $\gamma$ and $c_\epsilon$, for each $\epsilon \in
\{0,1\}^m$, $|\epsilon| \geq 0$, such that
\be
\det(M - \lambda \II) = (-\lambda)^{k-1} \left( (- \gamma - \lambda) (
  \balpha - \blambda)^{\bi} - \sum_{|\epsilon| \geq 0} c_{\epsilon}
  (\balpha - \blambda)^{\bi - \epsilon} \right)
 .
\ee
\end{proposition}

\begin{remark}
The benefit of Proposition \ref{poly-prop} is that we can immediately
conclude that the spectrum of $M$ contains $k-1$ zeros,
and if $\lambda <0$, $(1/\beta_j - \lambda)>0$ implies
\be
(-\lambda)^{1-k}\det(M - \lambda \II) = 0
\ee
can be rearranged into an equation of the form
\be
f(\lambda ) = g( \lambda),
\ee
with $f(\lambda) = \lambda$ and
\be
g( \lambda)  = - \sum_{|\epsilon| \geq 0} c_\epsilon(\balpha -
\blambda)^{-\epsilon} - \gamma.
\ee
Now $f'(\lambda) = 1$ and $g'( \lambda) <0$ for $\lambda < 0$ implies $f$ and $g$ can
intersect at most at one point for $\lambda <0$.  But since the trace of $M$ is negative by
construction, we conclude there is precisely one negative eigenvalue.
\end{remark}

\begin{proof}[Proof of Proposition \ref{poly-prop}]

The choices of $\alpha_2, \ldots ,
\alpha_{m+1}$ depending on the other $\alpha_j$ were made precisely so
that through row operations we can reduce $M$ to a matrix
\be
PM = \left( \begin{array}{cc} -\alpha \11_{1 \times k}  & \00_{1
      \times i} \\
    \00_{(k-1) \times k } & \00_{(k-1) \times i} \\
    M_{2,1}' & M_{2,2}' \end{array} \right),
\ee
where $(M_{2,1}', M_{2,2}')$ is the unchanged last $i$ rows from $M$.
Here $P$ is the invertible matrix whose action by left multiplication
is these row operations.  Computing the action of $P^{-1}$ by right
multiplication produces the corresponding similarity transformation,
and $M$ is similar to a matrix
\be
P M P^{-1} = \left( \begin{array}{ccc}
    -k \alpha & - \alpha \11_{1 \times (k-1)} & r_{1,3} \\
    \00_{(k-1) \times 1} & \00_{(k-1) \times ( k-1)} & \00_{(k-1)
      \times i} \\ c_{1,3} & M_{3,2}'' & \tB \end{array} \right),
\ee
where
\ben
r_{1,3}  & = & \left( \beta_1 \alpha \11_{1 \times i_1}, \beta_2 \alpha
  \11_{1 \times i_2}, \ldots , \beta_m \alpha \11_{1 \times i_m}
\right), \nonumber \\
c_{1,3} & = & - \11_{i \times 1}, \label{c13}
\een
and $M_{3,2}''$ is the $i \times (k-1)$ sub-matrix remaining unchanged
from $M$.  Here, \eqref{c13} follows from \eqref{choose-alpha} and $\tB$
is the $i \times i$ matrix given in row blocks:
\ben
\label{B-def}
\tB = \left( \begin{array}{c} \tR_1 \\ \tR_2 \\ \vdots \\ \tR_m \end{array}
\right),
\een
with
\be
\tR_1 & = & \left( \11_{i_1 \times i_1} + \beta_1^{-1} \II_{i_1 \times
    i_1}, \frac{\beta_2}{\beta_1} \11_{i_1 \times i_2} ,
  \frac{\beta_3}{\beta_1} \11_{i_1 \times i_3}, \ldots ,
  \frac{\beta_m}{\beta_1} \11_{i_1 \times i_m} \right) \\
\tR_2 & = & \left( \11_{i_2 \times i_1} , \11_{i_2 \times i_2} +
  \beta_2^{-1} \II_{i_2 \times i_2}, \frac{\beta_3}{\beta_2} \11_{i_2
    \times i_3} , \ldots, \frac{ \beta_m}{\beta_2} \11_{i_1 \times i_m} \right), \\
\tR_3 & = & \left( \11_{i_3 \times i_1}, \11_{i_3 \times i_2}, \11_{i_3
    \times i_3} + \beta_3^{-1} \II_{i_3 \times i_3} ,
  \frac{\beta_4}{\beta_3} \11_{i_3 \times i_4}, \ldots,
  \frac{\beta_m}{\beta_3} \11_{i_3 \times i_m} \right), \\
& \vdots & \\
\tR_{m-1} & = & \Bigg( \11_{i_{m-1} \times i_{1}}, \ldots , \11_{i_{m-1}
    \times i_{m-2}},\\
&& \quad \11_{i_{m-1} \times i_{m-1}} + \beta_{m-1}^{-1}
    \II_{i_{m-1} \times i_{m-1}}, \frac{\beta_m}{\beta_{m-1}}
    \11_{i_{m-1} \times i_m} \Bigg), \\
\tR_m & = & \left( \11_{i_m \times i_1}, \ldots, \11_{i_m \times
      i_{m-1}}, \11_{i_m\times i_m} + \beta_m^{-1} \II_{i_m \times
      i_m} \right).
\ee

Now $P M P^{-1}$ has rows of zero for $r_2$ through $r_k$.  Since
similarity transformations leave the spectrum invariant, we have
\be
\det(M - \lambda \II) & = & \det(PMP^{-1} - \lambda \II) \\
& = & (- \lambda)^{k-1}\det(M_0),
\ee
where $M_0$ is the $(i+1) \times (i+1)$ matrix
\be
M_0 = \left( \begin{array}{cc} -k \alpha - \lambda & r_{1,3} \\
    c_{1,3} & \tB- \lambda \II_{i \times i} \end{array} \right).
\ee
The following lemma is the induction step of the proof.  In order to
simplify notation, set
\ben
\label{B-def-2}
B = \tB - \11_{i \times i},
\een
and for $tR_j$ defined above,
\ben
\label{Rj-def}
R_j = \tR_j - \11_{i_j \times i}.
\een

\begin{lemma}
\label{M0-lemma}
For each $1 \leq j \leq m-1$, let $I_j = i_m + i_{m-1} + \ldots +
i_{m-j+1}$ and $\bI_j = (i_m, i_{m-1}, \ldots, i_{m-j+1})$.  There exist positive constants
$c_\epsilon$ and $(i - I_j +1) \times (i - I_j + 1)$ matrices
$M_\epsilon^j$ for each $\epsilon \in \{0,1\}^j$,
$c_0 = 1$,
 such that
\ben
\label{M-poly-def}
\det (M_0) = \sum_{\epsilon \in \{0,1\}^j} c_\epsilon (
\balpha - \blambda)^{(\bI_j - \epsilon)} \det(M_\epsilon^j).
\een

The $M_\epsilon^j$ satisfy the following properties:

(i)
\ben
\label{M0j-def}
 M_{0}^j = \left( \begin{array}{cc}
    -\gamma_j- \lambda & r_j \\
    -\11_{(i-I_j)  \times 1} & B^j - \lambda \II_{(i-I_j) \times
      (i-I_j)} \end{array} \right),
\een
where
\ben
\label{gamma-j}
\gamma_j = k \alpha - I_j >0,
\een
\be
r_j & = & \big( (\alpha(\beta_1-k) - \lambda)  \11_{1 \times i_1}, (\alpha(\beta_2-k) - \lambda) \11_{1
    \times i_2}, \\
&& \quad \ldots, (\alpha(\beta_{m-j}-k) - \lambda) \11_{1 \times i_{m-j}} \big),
\ee
and $B^j$ is the $(i-I_j) \times (i-I_j)$ matrix obtained from $B$ in \eqref{B-def-2} by removing
the last $I_j$ rows and columns.

(ii) For $\epsilon \neq 0$,
\ben
\label{Meps-def}
M_\epsilon^j = \left( \begin{array}{cc} - \gamma^j_\epsilon & r_j \\
    c^j_\epsilon & B^j - \lambda \II_{(i-I_j) \times (i-I_j)} \end{array} \right),
\een
where $r_j$ and $B_j$ are as in (i), $\gamma_\epsilon^j >0$, and
\be
c^j_\epsilon = \left( \begin{array}{c} - \delta_\epsilon^1 \11_{i_1
      \times 1} \\ - \delta_\epsilon^2 \11_{i_2 \times 1} \\ \vdots
    \\- \delta_\epsilon^{m-j} \11_{i_{m-j} \times 1} \end{array}
\right),
\ee
for constants
\be
\delta_\epsilon^p >0.
\ee

(iii) We have the relations
\be
\delta_{(\epsilon,1)}^p & = & \left( 1 - \frac{\beta_{m-j}}{\beta_p}
\right),  \\
\delta_{(\epsilon,0)}^p & = & \delta_\epsilon^p,
\ee
for $1 \leq p \leq m-j-1$, and
\be
-\gamma^{j+1}_{(\epsilon,1)} & = & \alpha (\beta_{m-j} - k) -
\beta_{m-j}^{-1}, \\
-\gamma^{j+1}_{(\epsilon,0)} & = & -\gamma_{\epsilon}^j + i_{m-j} \delta_\epsilon^{m-j}.
\ee
\end{lemma}
\begin{proof}
The basic idea is to inductively use row and column operations on the
last remaining row block.  At each step this results in the last row
block being diagonal except for negative elements in the
first column.  When we expand the determinant along the last row
block, we get terms involving $(\beta_j^{-1} - \lambda)^{i_{j}}$ and
terms involving $(\beta_j^{-1} - \lambda)^{i_j -1}$.  At each step in the induction, we then
permute the last column to the first, keeping track of the powers of
$-1$, to get a leading element of the form $-\gamma_j''$.  We will
carefully define all of this in the remainder of the proof.

{\bf Base case.}  The last row block of $M_0$ is of the form
\be
\left( - \11_{i_m \times 1}, R_m^0 \right),
\ee
where
\be
R_m^0  =  \left( \11_{i_m \times i_1}, \ldots, \11_{i_m \times
      i_{m-1}}, \11_{i_m\times i_m} + \beta_m^{-1} \II_{i_m \times
      i_m} \right).
\ee
To simplify this expression, we add column $1$ to the last $i$
columns, which has the effect of replacing $\tB$ with $B$ as defined
in \eqref{B-def-2}, but adds $-k \alpha - \lambda$ to the last $i$ elements in
the first row.  To eliminate the $-\lambda$ in the last $i_m$ elements
in the first row, we subtract the last
$i_m$ rows from the first and obtain
\be
\det(M_0) = \det(M_0'),
\ee
where
\be
M_0' = \left( \begin{array}{cc} -\gamma_1 - \lambda & r_1' \\
-\11_{i \times 1} & B \end{array} \right).
\ee
Here
\be
- \gamma_1  =  -k \alpha + i_m,
\ee
\be
r_1' & = & \big( (\alpha(\beta_1-k) - \lambda)  \11_{1 \times i_1}, (\alpha(\beta_2-k) - \lambda) \11_{1
    \times i_2}, \\
&& \quad \ldots, (\alpha(\beta_{m-1}-k) - \lambda) \11_{1 \times i_{m-1}} , - \gamma_1' \11_{1 \times i_m}
\big),
\ee
with
\be
-\gamma_1' = \alpha (\beta_m - k)  - \beta_m^{-1},
\ee
and $B$ is as defined in \eqref{B-def-2}.

Now when we expand the determinant of $M_0'$ along the last row block,
we get contributions from the $-1$s in the first column, and the
diagonal elements in $R_m$.  Expanding the determinant along these
rows and permuting the resulting submatrices so that the lower right
$(i-i_m) \times (i-i_m)$ submatrix agrees with the definition of
$B^1$ in the lemma, we obtain
\be
\det(M_0') = (\beta_m^{-1} - \lambda)^{i_m} \det ( M_1 ) + i_m
(\beta_m^{-1} - \lambda)^{i_m -1} \det(M_1'),
\ee
where
\be
M_1 = \left( \begin{array}{cc} - \gamma_1 - \lambda & r_1 \\
    - \11_{(i-i_m) \times 1} & B^1 \end{array} \right)
\ee
and
\be
M_1' = \left( \begin{array}{cc} - \gamma_1' & r_1 \\ c_1 &
    B^1 \end{array} \right).
\ee
Here $B^1$ is defined in the statement of the lemma,
\be
r_1 = \left( \beta_1 \alpha, \beta_2 \alpha, \ldots, \beta_{m-1}
  \alpha \right),
\ee
and
\be
c_1 = \left( \begin{array}{c} - \delta_1^1 \11_{i_1 \times 1} \\ -
    \delta_2^1 \11_{i_2 \times 1} \\ \vdots \\ -\delta_{m-1}^1
    \11_{i_{m-1} \times 1} \end{array} \right),
\ee
with
\be
\delta_p^1 = 1-\frac{\beta_m}{\beta_p} ,
\ee
in accordance with the statement of the Lemma.

{\bf Induction step.}  Now suppose the Lemma is true for some $1 \leq
j \leq m-2$.  We show the same reduction ideas used in the base case
will reduce to the statement of the Lemma for $j+1$.  That is, assume
we have matrices and constants as in the formula \eqref{M-poly-def}.
For $M_0^j$, we subtract the
last $i_{m-j}$ rows from the first to eliminate the $-\lambda$ in the
last $i_{m-j}$ elements in the first row.  Expanding the determinant
along the last row block (the last $i_{m-j}$ rows) and permuting as
necessary yields
\be
\det (M_0^j) = ( \beta_{m-j}^{-1} - \lambda )^{i_{m-j}} \det (M_0^{j+1}) +
i_{m-j} ( \beta_{m-j}^{-1} - \lambda )^{i_{m-j} -1} \det ( M_1^{j+1}
),
\ee
where $M_0^{j+1}$ is defined in the statement of the Lemma, and
\be
M_1^{j+1} = \left(\begin{array}{cc} - \gamma_{j+1}^1 & r_{j+1} \\ c_{j+1}^1
    & B^{j+1} \end{array} \right),
\ee
where
\be
- \gamma_{j+1}^1 = \alpha(\beta_{m-j} -k) - \beta_{m-j}^{-1},
\ee
$r_{j+1}$ is as defined in the Lemma, and
\be
c_{j+1}^1 = \left( \begin{array}{c} - \delta_1^{j+1,1} \11_{i_1 \times
      1} \\ - \delta_2^{j+1,1} \11_{i_2 \times 1} \\ \vdots \\
    -\delta_{m-j-1}^{j+1,1} \11_{i_{m-j-1} \times 1} \end{array}
\right),
\ee
with
\ben
\label{delta-pj1}
-\delta_{p}^{j+1,1} = \frac{\beta_{m-j}}{\beta_p} -1 <0.
\een
Observe the multi-index $\epsilon_1$ associated to $M_1^{j+1}$ is
\be
\epsilon_1 = (\underbrace{ 0, \ldots, 0 }_{j}, 1),
\ee
so $\delta_{\epsilon_1}^p := \delta_p^{j+1,1}$ in \eqref{delta-pj1}
agrees with the statement of the Lemma.

We next tackle $M_\epsilon^j$ for $|\epsilon|>0$.  We again subtract
the last $i_{m-j}$ rows from the first row to eliminate the
$-\lambda$s in the last $i_{m-j}$ columns.  As with the other cases, we
then expand the determinant along the last row block and permute as
necessary to ensure $B^{j+1}$ be the lower right submatrix.  We have
\be
\lefteqn{\det{M_\epsilon^j}= } \\
&&  = (\beta_{m-j}^{-1} - \lambda )^{i_{m-j}}
\det(M_{\epsilon_2}^{j+1}) + i_{m-j} \delta_\epsilon^{m-j} (
\beta_{m-j}^{-1} - \lambda)^{i_{m-j} -1} \det(M_{\epsilon_3}^{j+1}),
\ee
where
\be
M_{\epsilon_2}^{j+1} = \left( \begin{array}{cc} - \gamma_\epsilon^j &
    r_{j+1} \\
    c_{\epsilon_2}^{j+1} & B^{j+1} \end{array} \right),
\ee
with $\gamma_\epsilon^j$ inherited from $M_\epsilon^j$, $r_{j+1}$ as
defined in the Lemma, and
\be
c_{\epsilon_2}^{j+1} = \left( \begin{array}{c}  - \delta_\epsilon^1 \11_{i_1
      \times 1} \\ - \delta_\epsilon^2 \11_{i_2 \times 1} \\ \vdots
    \\- \delta_\epsilon^{m-j-1} \11_{i_{m-j-1} \times 1} \end{array}
\right),
\ee
with $\delta_\epsilon^p$ inherited from $M_\epsilon^j$.  Here
\be
M_{\epsilon_3}^{j+1} = \left( \begin{array}{cc} -
    \gamma_{\epsilon_3}^{j+1} & r_{j+1} \\
    c_{\epsilon_3}^{j+1} & B^{j+1} \end{array} \right),
\ee
where
\be
-\gamma_{\epsilon_3}^{j+1} = \alpha(\beta_{m-j} - k) - \beta_{m-j}^{-1} ,
\ee
and
\be
c_{\epsilon_3}^{j+1} = \left( \begin{array}{c}  - \delta_{\epsilon_3}^{j+1,1} \11_{i_1
      \times 1} \\ - \delta_{\epsilon_3}^{j+1,2} \11_{i_2 \times 1} \\ \vdots
    \\- \delta_{\epsilon_3}^{j+1, m-j-1} \11_{i_{m-j-1} \times 1} \end{array}
\right),
\ee
with
\be
-\delta_{\epsilon_3}^{j+1,p} = \frac{\beta_{m-j}}{\beta_{p}} - 1,
\ee
as in the statement of the Lemma.

Relabeling as necessary, this completes the proof of the induction step, and hence the proof of
the Lemma.
\end{proof}

Now in order to finish the proof of Proposition \ref{poly-prop} we
calculate what happens in the $j = m$ case.  From Lemma
\ref{M0-lemma}, we have
\be
\det (M_0) = \sum_{\epsilon \in \{0,1\}^{m-1}} c_\epsilon (
\balpha - \blambda)^{(\bI_{m-1} - \epsilon)} \det(M_\epsilon^{m-1}),
\ee
where, relabeling for simplicity in exposition,
\be
M_0^{m-1} = \left( \begin{array}{cc} - \gamma - \lambda &
    (\alpha(\beta_1 -k) - \lambda) \11_{1 \times i_1} \\
    -\11_{i_1 \times 1} &  (\beta_1^{-1} -
    \lambda) \II_{i_1 \times i_1} \end{array} \right),
\ee
with $-\gamma = -  k\alpha + i - i_1$, and
\be
M_\epsilon^{m-1} = \left( \begin{array}{cc} -\gamma_\epsilon &  (\alpha(\beta_1 -k) - \lambda)  \11_{1 \times i_1} \\ - \delta_\epsilon \11_{i_1 \times 1}
    &  (\beta_1^{-1} -
    \lambda) \II_{i_1 \times i_1} \end{array} \right).
\ee

Using {\it iii} from the Lemma, if $\epsilon = (\epsilon',1)$ with
$\epsilon' \in \{0,1\}^{m-2}$,
\be
-\gamma_{(\epsilon',1)} = \alpha(\beta_2 - k) - \beta_2^{-1},
\ee
and if $\epsilon = (\epsilon',0)$,
\be
-\gamma_{(\epsilon',0)} & = & -\gamma_{(\epsilon',0)}^{m-1} \\
& = & - \gamma_{\epsilon'}^{2} + i_2 \delta_{\epsilon'}^2 \\
& < & -\gamma_{\epsilon'}^{2} + i_2 \\
& \leq & -\gamma^j,
\ee
by induction.


Proceeding as in the proof of the Lemma
yields
\be
\det (M_0^{m-1}) = (- \gamma' - \lambda ) (\beta_1^{-1} - \lambda
)^{i_1} -i_1 \gamma'' (\beta_1^{-1} - \lambda )^{i_1 -1},
\ee
where
\be
- \gamma' & = & -k \alpha + i \text{ and} \\
- \gamma'' & = & \alpha(\beta_1 -k) - \beta_1^{-1}
\ee
both of which are negative.  Similarly,
\be
\det ( M_\epsilon^{m-1}) = - \gamma_{\epsilon} (\beta_1^{-1} -
\lambda )^{i_1} - i_m \delta_ \epsilon \gamma_\epsilon' (
\beta_1^{-1} - \lambda)^{i_1 -1},
\ee
where $\gamma_\epsilon'$ was defined above, and is negative.  This
proves the Proposition.

\end{proof}

\begin{remark}
We illustrate the proof of Proposition \ref{poly-prop} by following
the steps in the concrete example of $M = $\eqref{cdv-2}.  We have
\be
\det(M - \lambda \II) = \det(PMP^{-1} - \lambda \II),
\ee
with
\be
\lefteqn{PMP^{-1} =} \\ && \left( \begin{array}{cccccccccccc}
    -7 \alpha & - \alpha & -\alpha & - \alpha & - \alpha & - \alpha &
    - \alpha & 5 \alpha & 5 \alpha & 4 \alpha & \alpha & \alpha \\
    0 & 0 & 0 & 0 & 0 & 0 & 0 & 0 & 0 & 0 & 0 & 0 \\
    0 & 0 & 0 & 0 & 0 & 0 & 0 & 0 & 0 & 0 & 0 & 0 \\
    0 & 0 & 0 & 0 & 0 & 0 & 0 & 0 & 0 & 0 & 0 & 0 \\
    0 & 0 & 0 & 0 & 0 & 0 & 0 & 0 & 0 & 0 & 0 & 0 \\
    0 & 0 & 0 & 0 & 0 & 0 & 0 & 0 & 0 & 0 & 0 & 0 \\
    0 & 0 & 0 & 0 & 0 & 0 & 0 & 0 & 0 & 0 & 0 & 0 \\
-1 & -\frac{1}{5} & -\frac{1}{5} & -\frac{1}{5} & -\frac{1}{5} & 0 & 0 & \frac{6}{5} & 1 & \frac{4}{5} & \frac{1}{5} & \frac{1}{5} \\
-1 & -\frac{1}{5} & -\frac{1}{5} & -\frac{1}{5} & -\frac{1}{5} & 0 & 0
& 1 & \frac{6}{5} & \frac{4}{5} & \frac{1}{5} & \frac{1}{5} \\
-1 & -\frac{1}{4} & -\frac{1}{4} & -\frac{1}{4} & 0 & 0 & 0 & 1 & 1 & \frac{5}{4} & \frac{1}{4} & \frac{1}{4} \\
-1 & 0 & 0 & 0 & 0 & 0 & 0 & 1 & 1 & 1 & 2 & 1 \\
-1 & 0 & 0 & 0 & 0 & 0 & 0 & 1 & 1 & 1 & 1 & 2
\end{array} \right),
\ee
where $\alpha = \alpha_1 + 53/20$ for $\alpha_1>0$ yet to be
determined.  The rows of zeros yield a contribution to $\det(M -
\lambda \II)$ of
$(- \lambda)^6$, and we have reduced to studying $\det(M_0)$, for
\be
M_0 = \left( \begin{array}{cccccc}
    -7 \alpha - \lambda & 5 \alpha &  5 \alpha & 4 \alpha & \alpha &
    \alpha \\
    -1 & \frac{6}{5}- \lambda & 1 & \frac{4}{5} & \frac{1}{5} & \frac{1}{5} \\
    -1 & 1 & \frac{6}{5}- \lambda  & \frac{4}{5} & \frac{1}{5} & \frac{1}{5}\\
    -1 & 1 & 1 & \frac{5}{4} - \lambda & \frac{1}{4} & \frac{1}{4} \\
    -1 & 1 & 1 & 1& 2-\lambda & 1 \\
    -1 & 1 & 1 & 1 & 1 & 2 - \lambda \end{array} \right).
\ee
We use the first column to kill the $1$s in the last three rows, and
then the last two rows to kill the $- \lambda$s in the first row to
get $\det(M_0) = \det(M_0')$, where
\be
\lefteqn{ M_0' = } \\ && \left( \begin{array}{cccccc}
    -7 \alpha +2 - \lambda & -2 \alpha - \lambda &  -2 \alpha -
    \lambda & -3 \alpha - \lambda & -6\alpha-1 &
    -6\alpha -1 \\
    -1 & \frac{1}{5}- \lambda & 0 & -\frac{1}{5} & -\frac{4}{5} & -\frac{4}{5} \\
    -1 & 0 & \frac{1}{5}- \lambda  & -\frac{1}{5} & -\frac{4}{5} & -\frac{4}{5}\\
    -1 & 0 & 0 & \frac{1}{4} - \lambda & -\frac{3}{4} & -\frac{3}{4} \\
    -1 & 0 & 0 & 0 & 1-\lambda & 0 \\
    -1 & 0 & 0 & 0 & 0 & 1 - \lambda \end{array} \right).
\ee
Expanding $\det(M_0')$ along the last two lines yields:
\be
\det(M_0') =  2 (1 - \lambda) \det(M_1) + (1-\lambda)^2 \det(M_2),
\ee
where
\be
M_1 = \left( \begin{array}{cccc}
-6 \alpha -1 & -2 \alpha - \lambda & -2 \alpha - \lambda & -3 \alpha -
\lambda \\
-\frac{4}{5} & \frac{1}{5}- \lambda & 0 & -\frac{1}{5} \\
-\frac{4}{5} & 0 & \frac{1}{5}- \lambda  & -\frac{1}{5} \\
- \frac{3}{4} & 0 & 0 & \frac{1}{4} - \lambda \end{array} \right),
\ee
and
\be
M_2 = \left( \begin{array}{cccc}
-7 \alpha + 2 - \lambda & -2 \alpha - \lambda & -2 \alpha - \lambda &
-3 \alpha  - \lambda \\
-1 & \frac{1}{5}- \lambda & 0 & -\frac{1}{5} \\
-1 & 0 & \frac{1}{5}- \lambda  & -\frac{1}{5} \\
- 1 & 0 & 0 & \frac{1}{4} - \lambda \end{array} \right).
\ee
Subtract the last row
from the first in $M_1$ and $M_2$ and expanding the respective determinants along the last
row yields
\be
\det(M_1) & = & \frac{3}{4} \det(M_3) + \left( \frac{1}{4} - \lambda \right)
  \det (M_4), \\
\det(M_2) & = & \det(M_5) + \left( \frac{1}{4} - \lambda \right) \det (M_6),
\ee
where
\be
M_3 & = & \left( \begin{array}{ccc}
-3 \alpha - \frac{1}{4} & -2 \alpha -\lambda & -2 \alpha -\lambda \\
-\frac{1}{5} & \frac{1}{5}- \lambda & 0  \\
-\frac{1}{5} & 0 & \frac{1}{5}- \lambda
\end{array} \right), \\
M_4 & = & \left( \begin{array}{ccc}
-6 \alpha - 1 & -2 \alpha -\lambda & -2 \alpha -\lambda \\
-\frac{4}{5} & \frac{1}{5}- \lambda & 0  \\
-\frac{4}{5} & 0 & \frac{1}{5}- \lambda
\end{array} \right), \\
M_5 & = & \left( \begin{array}{ccc}
-3 \alpha - \frac{1}{4} & -2 \alpha -\lambda & -2 \alpha -\lambda \\
-\frac{1}{5} & \frac{1}{5}- \lambda & 0  \\
-\frac{1}{5} & 0 & \frac{1}{5}- \lambda
\end{array} \right), \text{ and }\\
M_6 & = & \left( \begin{array}{ccc}
-7 \alpha +3 - \lambda & -2 \alpha - \lambda & -2 \alpha - \lambda \\
-1 & \frac{1}{5}- \lambda & 0  \\
-1 & 0 & \frac{1}{5}- \lambda
\end{array} \right).
\ee
Expanding the determinant as in the previous iterations yields:
\be
\det(M_3) & = & \frac{2}{5} (-2 \alpha - \frac{1}{5}) \left( \frac{1}{5} - \lambda
\right) + \left(- 3 \alpha + \frac{3}{20} \right) \left( \frac{1}{5} - \lambda
\right)^{2}, \\
\det(M_4) & = & \frac{8}{5} \left(-2 \alpha - \frac{1}{5} \right) \left( \frac{1}{5} - \lambda
\right) + \left(- 6 \alpha + \frac{27}{20} \right) \left( \frac{1}{5} - \lambda
\right)^{2}, \\
\det(M_5) & = & \det(M_3), \text{ and} \\
\det(M_6) & = & 2 \left(-2 \alpha - \frac{1}{5} \right) \left( \frac{1}{5} - \lambda
\right) + \left(- 7 \alpha +5 -\lambda \right) \left( \frac{1}{5} - \lambda
\right)^{2}.
\ee
Following our calculations back to the original matrix $M$, we see
$\det(M)$ satisfies Proposition \ref{poly-prop} as long as $\alpha_1>0$
is chosen large enough so that $\alpha = \alpha_1 + 53/20$ satisfies
the inequalities
\be
-7 \alpha + 5 & < & 0, \\
-3 \alpha + \frac{3}{20} & < & 0, \text{ and}\\
-6 \alpha + \frac{27}{20} & < & 0.
\ee


\end{remark}
\qed

{\bf Proof of property \eqref{M3}.}

In order to verify \eqref{M3}, we write such a matrix $X$ as
\be
X = \left(  \tC_1, \ldots, \tC_{2m+1} \right)
\ee
with $\tC_j$ having the same dimensions as $C_j$ from $M$.  We also
let $\tR_j = \tC_j^T$ be the row blocks of $X$.  The assumptions on
$X$ and the construction of $M$ show
\be
\tC_1 = \00_{n \times k_{m+1}}
\ee
and
\be
\tC_2 = \left( \begin{array}{c} \00_{(n-i_m) \times k_m} \\ \tC^2_{i_m
      \times k_m} \end{array} \right),
\ee
where $\tC^2_{i_m \times k_m}$ denotes an arbitary $i_m \times k_m$
matrix
particular to $\tC_2$.
The equation $MX=0$ implies
\be
R_{2m+1} \tC_2 = \00_{i_m \times k_m},
\ee
which further implies
\be
 \II_{i_m \times i_m} \tC^2_{i_m \times k_m} = \00_{i_m
  \times k_m},
\ee
so $\tC_2$ is zero, as well as $\tR_2 = \tC_2^T$.  Continuing like this, for $3
\leq j \leq m+1$
\be
R_{2m+1} \tC_j = \00_{i_m \times k_{m+2 - j}},
\ee
and for $m+2 \leq j \leq 2m+1$,
\be
R_{2m+1} \tC_j = \00_{i_m \times i_{j-m-1}},
\ee
together imply $\tR_{2m+1} = \00_{i_m \times n}$.  Now for the
purposes of induction, suppose we know $\tR_{2m+2 - j} = \00_{i_{m+2-j}
  \times n}$, $\tC_{j+1} = \00_{n \times k_{m+1-j}}$ for some $2 \leq
j \leq m-1$.  Multiplying $\tC_{j+1+k}$ on the left by $R_{2m+1-j}$ and
setting equal to zero for $1 \leq k \leq m-j$ gives $\tR_{2m+1-j} =
\00_{i_{m-j} \times n}$ and hence $\tC_{j+2} = \00_{n \times
  k_{m-j}}$.  Thus by induction $X = 0$ and $M$ satisfies \eqref{M3}.

\begin{remark}
In order to illustrate the proof of property \eqref{M3} for $M$, we
show how it works for our example \eqref{cdv-2}.  For this $M$, $X$
has the form
\be
X = \left( \begin{array}{rrrrrrrrrrrr}
0 & 0 & 0 & 0 & 0 & 0 & 0 & 0 & 0 & 0 & 0 & 0 \\
0 & 0 & 0 & 0 & 0 & 0 & 0 & 0 & 0 & 0 & x_1 & x_2 \\
0 & 0 & 0 & 0 & 0 & 0 & 0 & 0 & 0 & 0 & x_3 & x_4 \\
0 & 0 & 0 & 0 & 0 & 0 & 0 & 0 & 0 & 0 & x_5 & x_6 \\
0 & 0 & 0 & 0 & 0 & 0 & 0 & 0 & 0 & x_7 & x_8 & x_9 \\
0 & 0 & 0 & 0 & 0 & 0 & 0 & x_{10} & x_{11} & x_{12} & x_{13} & x_{14} \\
0 & 0 & 0 & 0 & 0 & 0 & 0 & x_{15} & x_{16} & x_{17} & x_{18} & x_{19}
\\
0 & 0 & 0 & 0 & 0 & x_{10} & x_{15} & 0 & x_{20} & x_{21} & x_{22} &
x_{23} \\
0 & 0 & 0 & 0 & 0 & x_{11} & x_{16} & x_{20} & 0 & x_{24} & x_{25} &
x_{26} \\
0 & 0 & 0 & 0 & x_7 & x_{12} & x_{17} & x_{21} & x_{24} & 0 & x_{27} &
x_{28} \\
0 & x_1 & x_3 & x_5 & x_8 & x_{13} & x_{18} & x_{22} & x_{25} & x_{27}
& 0 & x_{29} \\
0 & x_2 & x_4 & x_6 & x_9 & x_{14} & x_{19} & x_{23} & x_{26} & x_{28}
& x_{29} & 0
\end{array} \right).
\ee
Multiplying $X$ on the left by the last two rows of $M$ implies the
last two rows of $X$ are zero.  Hence the last two columns of $X$ are
zero, which from the structure of $X$ implies the first $4$ rows of
$X$ are zero.  Now the $10$th row of $M$ has four non-zero entries
followed zeros and a non-zero entry in the tenth position.  This
implies row $10$ of $X$ is zero, so column $10$ is zero, and hence row
$5$ is zero.  Continuing in this fashion eventually gives $X=0$.
\end{remark}

\qed

Following the numbering schemes used in the proof of Proposition
\ref{poly-prop} we get the following estimates on the size of
the eigenvalues of $M$.
\begin{corollary}
If $\lambda<0$ is the negative eigenvalue of $M$, then
\be
\lambda < -k \alpha + i.
\ee
If $\lambda >0$ is a positive eigenvalue of $M$, then
\be
\lambda \geq \beta_1^{-1} > \frac{1}{k}.
\ee
\end{corollary}

\end{document}